\documentclass[11pt, a4paper]{article}
\usepackage{amsmath,amssymb,amsthm,bm, mathrsfs}
\usepackage{subcaption}
\usepackage[affil-it]{authblk}
\usepackage{hyperref,graphicx}
\usepackage[margin=2.6cm]{geometry}

\newcommand{\bu}{\bm{u}}
\newcommand{\bv}{\bm{v}}
\newcommand{\bnu}{\bm{\nu}}
\newcommand{\bpi}{\bm{\pi}}
\newcommand{\bx}{\bm{x}}

\newcommand{\demi}{\frac{1}{2}}
\begin{document}
\title{Lagrange-flux schemes: reformulating second-order accurate Lagrange-remap schemes for better node-based HPC performance}
%
% if title is too large for page header
%
%%%\titlerunning{Short title}
%
% subtitle is optionnal
%
%%%\subtitle{Do you have a subtitle?\\ If so, write it here}
%\author{Florian De Vuyst\inst{\ref{un}}\correspondingauthor \and
        %% corresponding author following by \correspondingauthor
        %Thibault Gasc\inst{\ref{un},\ref{deux},\ref{trois}} \and
				%Renaud Motte\inst{\ref{trois}} \and
				%Mathieu Peybernes\inst{\ref{quatre}} \and
				%Rapha\"el Poncet\inst{\ref{cinq}}
        %% etc.
%}
\author[1]{Florian De Vuyst \thanks{Electronic address: devuyst@cmla.ens-cachan.fr}}
%
%\thanks{}
\affil[1]{Center for Mathematics and their applications CMLA, \'Ecole Normale
Sup\'erieure de Cachan, Universit\'e Paris-Saclay, CNRS, 94235 Cachan France}

\author[1,2,3]{Thibault Gasc}
\affil[2]{Maison de la Simulation, URS 3441, CEA Saclay, 91191 Gif-sur-Yvette, France}

\author{Renaud Motte}
\affil[3]{CEA DAM DIF, F-91297 Arpajon France}

\author{Mathieu Peybernes}
\affil[4]{CEA Saclay, DEN, 91191 Gif-sur-Yvette France}

\author{Raphael Poncet}
\affil[5]{CGG, 27 avenue Carnot, 91300 Massy France}
%
% if author is too large for page header
%
%\authorrunning{F. De Vuyst et al.}
%%
%\institute{CMLA, ENS Cachan, Universit\'e Paris-Saclay, CNRS, 94235 Cachan France\email{devuyst@cmla.ens-cachan.fr}\label{un} \and
           %\email{thibault.gasc@cea.fr}\label{deux} \and
           %CEA DAM DIF, F-91297 Arpajon France\email{renaud.motte@cea.fr}\label{trois} \and
					 %CEA Saclay, DEN, 91191 Gif-sur-Yvette, France\email{mathieu.peybernes@cea.fr}\label{quatre} \and
					 %CGG, 27 avenue Carnot, 91300 Massy France\email{raphael.poncet@cgg.com}\label{cinq}
          %}
%
%\resume{%
  %%\frenchtitle{abstract title (optionnal)}
  %Dans un article récent \cite{PARCO2015}, nous avons effectué l'analyse de la performance d'un schéma de type
	%Lagrange+projection à variables décalées~; cette classe de solveurs est très utilisée pour les applications
	%d'hydrodynamique matériau. Dans cet article, on s'intéresse à la reformulation des solveurs Lagrange-projection
	%afin d'améliorer leur performance globale sur architectures de calcul standards. De manière inattendue,
	%l'analyse nous a conduit vers la découverte d'une nouvelle famille de solveurs -- appelés schémas
	%Lagrange-flux -- qui apparaissent comme très prometteurs dans la communauté CFD.
%}
\maketitle

\begin{abstract}%
  %\englishtitle{abstract title (optionnal)}
  In a recent paper \cite{PARCO2015}, we have achieved the performance analysis
of staggered Lagrange-remap schemes, a class of solvers widely used for Hydrodynamics applications. This paper is devoted to the rethinking and redesign of the Lagrange-remap process for achieving better performance using today's computing architectures.  As an unintended outcome, the analysis has lead us to the discovery of a new family
of solvers -- the so-called Lagrange-flux schemes -- that appear to be promising for the CFD community.
\end{abstract}

\section{Motivation and introduction}
\begin{figure}
\begin{center}
%\notoprule
\includegraphics[width=0.9\linewidth]{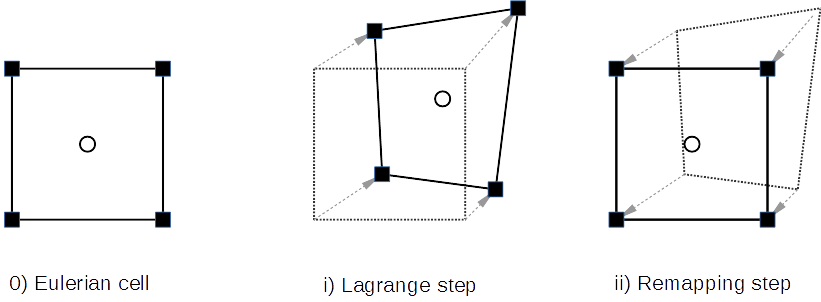}
\caption{``Legacy'' staggered Lagrange-remap scheme: thermodynamical variables are located at cell centers (circles) whereas velocity variables are located at cell nodes (squares).}\label{fig:1}
\end{center}
\end{figure}

For complex compressible flows involving multiphysics phenomenons like e.g.\ high-speed elastoplasticity, multimaterial interaction, plasma, gas-particles etc., 
a Lagrangian description of the flow is generally preferred. To achieve
robustness, some spatial remapping on a regular mesh may be added. A particular
case is the family of the so-called Lagrange+remap schemes \cite{Hirt1974,benson,youngs}, also referred to as remapped Lagrange solvers that apply a remap step on a reference (Eulerian) mesh after each Lagrangian time advance. Legacy codes
implementing remapped Lagrange solvers usually define thermodynamical variables
at cell centres and velocity at cell nodes (see figure~\ref{fig:1}).
In Poncet et al.~\cite{PARCO2015}, we have achieved a multicore node-based performance analysis of a reference Lagrange-remap hydrodynamics solver used in industry. By analyzing
each kernel of the whole algorithm, using roofline-type models \cite{Williams2009} on one side and
refined Execution Cache Memory (ECM) models  \cite{Treibig2010}, \cite{Stengel2015} on the other side, we have been able not only to quantitatively predict the performance of the whole algorithm --- with relative errors in the single digit range --- but also to identify a set of
features that limit the whole performance. This can be roughly summarized into
three points:
\begin{enumerate}
\item For typical mesh sizes of real applications, spatially staggered variables involve a rather big amount of communication to/from CPU caches and memory with low arithmetic intensity, thus lowering the whole performance;
\item Usual alternating direction (AD) strategies (see the appendix in~\cite{Collela1984}) or AD remapping procedures also generate too much communication with a loss of CPU occupancy.
\item For multimaterial flows using VOF-based interface reconstruction methods, 
there is a strong loss of performance due to some array indirections and noncoalescent
data in memory. Vectorization of such algorithms is also not trivial.\medskip
\end{enumerate}
From these observations and as a result of the analysis, we decided to ``rethink''
Lagrange-remap schemes, with possibly modifying some aspects of the solver in order
to improve node-based performance of the hydrocode solver. We have looked for
alternative formulations that reduce communication and improve both arithmetic intensity
and SIMD (Single Instruction -- Multiple Data) property of the
algorithm. In this paper, we describe the process of redesign of Lagrange-remap
schemes leading to higher performance solvers. Actually, this redesign methodology also
gave us ideas of innovative eulerian solvers. The emerging methods, named Lagrange-flux schemes, appear to be promising in the extended computational fluid dynamics community.

The paper is organized as follows. In section~\ref{sec:2}, we first formulate the requirements for the design of better Lagrange-remap schemes. In section~\ref{sec:3lag} we give a description of the Lagrange step and formulate it
under a finite volume form.
In section~\ref{sec:3}
we focus on the remap step which is reformulated as a finite volume scheme with
pure convective fluxes. This interpretation is applied in section~\ref{sec:4} to build
the so-called Lagrange-flux schemes. We also discuss the important issue of
achieving second order accuracy (in both space and time). In section~\ref{sec:5}
we comment the possible extension to multimaterial flow with the use of low-diffusive
and accurate interface-capturing methods. We will close the paper by some 
concluding remarks, work in progress and perspectives.
\section{Requirements}\label{sec:2}
%--------------------
%
Starting from ``legacy'' Lagrange-remap solvers and related observed performance
measurements, we want to improve the performance of these solvers by modifying some of the features of them but under some constraints and requirements:
\begin{enumerate}
\item A Lagrangian solver (or description) must be used (allowing for multiphysics coupling).
\item To reduce communication, we prefer using collocated cell-centered variables rather than
a staggered scheme.
\item To reduce communication, we prefer using a direct multidimensional remap solver rather than splitted alternating direction projections.
\item The method can be simply extended to second-order accuracy (in space and time).
\item The solver must be able to be naturally extended to multimaterial flows.  
\end{enumerate}
Before going further, let us comment the above requirements.
The second requirement should imply the use of a cell-centered Lagrange solver.
Fairly recently, Despr\'es and Mazeran in \cite{Despres2005} and 
Maire and et al.~\cite{Maire2007} (with high-order extension in~\cite{Maire2009}) 
have proposed pure cell-centered Lagrangian solvers based on 
the reconstruction of nodal velocities. In our study, we will examine if it is
possible to use approximate and simpler Lagrangian solvers in the Lagrange+remap context, in particular for the sake of performance. The fourth assertion
requires a full multidimensional remapping step, probably taking into
account geometric elements (deformation of cells and edges) if we want to ensure
high-order accuracy remapping. 
To summarize, our requirements are somewhat contradictory, and we have to find a good compromise between some simplifications-approximations and a loss of accuracy (or properties) of the numerical solver.
\section{Lagrange step} \label{sec:3lag}
%-----------------------------------
%
As example, let us consider the compressible Euler equations
for two-dimensional problems. Denoting $\rho,\ \bu=(u_i)_i,$ $i\in\{1,2\}$,\ $p$ and $E$ the density,
velocity, pressure and specific total energy respectively, the mass, momentum and
energy conservation equations read
\begin{equation}
\partial_t U_\ell + \nabla\cdot(\bu U_\ell) + \nabla\cdot\bpi_\ell=0,
\quad \ell=1,\dots,4,
\end{equation}
where $U=(\rho,(\rho u_i)_i,\rho E)$, $\bpi_1=\vec 0$, $\bpi_2=(p,0)^T$,
$\bpi_3=(0,p)^T$ and $\bpi_4=p\bu$.
For the sake of simplicity, we will use a perfect gas equation of state
$p=(\gamma-1)\rho (E-\frac{1}{2}|\bu|^2)$, $\gamma\in(1,3]$. 
The speed of sound $c$ is given by $c=\sqrt{\gamma p/\rho}$.

For any volume $V_t$ that moves with the fluid, from the Reynolds transport theorem
we have
\[
\frac{d}{dt} \int_{\mathscr{V}_t} U_\ell\, d\bx =
\int_{\partial \mathscr{V}_t} 
\left\{\partial_t U_\ell+\nabla\cdot (\bu\, U_\ell)\right\}\, d\bx 
= - \int_{\partial \mathscr{V}_t} \bpi_\ell\cdot\bnu \, d\sigma
\]
where $\bnu$ is the normal unit vector exterior to $\mathscr{V}_t$. This leads
to a natural explicit finite volume scheme in the form
\begin{eqnarray}
&&|K^{n+1,L}| (U_\ell)_K^{n+1,L} = |K| (U_\ell)_K^n \nonumber \\ [1.1ex]
&& \phantom{aaaaa} - \Delta t^n
\sum_{A^{n+\demi,L}\subset \partial K^{n+\demi,L}} |A^{n+\demi,L}|\, \bpi_A^{n+\demi,L}\cdot
\nu_A^{n+\demi,L}.
\label{eq:2}
\end{eqnarray}
In~\eqref{eq:2}, the superscript ``L'' indicates the Lagrange evolution of the quantity.
Any Eulerian cell $K$ is deformed into the Lagrangian volume $K^{n+\demi,L}$ at time
$t^{n+\demi}$, and into the Lagrangian volume $K^{n+1,L}$ at time $t^{n+1}$.
The pressure flux terms through the edges $A$ are evaluated at time $t^{n+\demi}$ in order to get second-order accuracy in time. Of course, that means that we need a predictor
step for the velocity field $\bu^{n+\demi,L}$ at time $t^{n+\demi}$ (not written here
for simplicity).

\paragraph{Notations.}  To simplify, in all what follows we will use the notation
$\bv^{n+\demi}=\bu^{n+\demi,L}$.
\section{Rethinking the remapping step}\label{sec:3}
%-----------------------------------------------
%
The remapping step considers the remapping the fields $U_\ell$ defined at cell centers
$K^{n+1,L}$ on the initial (reference) Eulerian mesh with cells $K$. 
Let us denote $\mathscr{R}^{n+1,L}$ a linear operator that reconstructs
piecewise polynomial functions from a discrete 
field $U^{n+1,L}$ defined at cell centers of the Lagrangian mesh 
$\mathscr{M}^{N+1,L}$ at time $t^{n+1}$.
The remapping process consists in projecting the reconstructed field on 
piecewise-constant function on the Eulerian mesh, according to the
intergral formula
\begin{equation}
U_K^{n+1} = \frac{1}{|K|}\int_{K}  \mathscr{R}^{n+1,L}U^{n+1,L}(\bx)\,d\bx.
\label{eq:3}
\end{equation}
Practically, they are many ways to consider the projection 
operation~\eqref{eq:3}. One can assemble elementary projection
contributions by computing the volume intersections between the 
reference mesh and the deformed mesh. But this procedure requires
the computation of all the geometrical elements. Moreover, the projection
needs local tests of projection with conditional branching (think about the 
very different cases of a compression and expansion). Thus the procedure 
is not SIMD and potentially leads to a loss of performance. 
The incremental remapping can also be interpreted as a transport/advection
algorithm, as emphasized by Dukowicz and Baumgardner~\cite{Dukowicz2000}
that appears to be better suited for SIMD treatments.

Let us now write a different original formulation of the remapping step. 
In this step, there is no time evolution of any quantity, and in some sense
we have $\partial_t U=0$, that we write
\[
\partial_t U = \partial_t U +\nabla\cdot(-\bv^{n+\demi} U)
\ +\ \nabla\cdot(\bv^{n+\demi} U) = 0.
\]
We decide to split up this equation into two substeps:
\begin{itemize}
\item [i)] ~Backward convection: 
\begin {equation}
\partial_t U+ \nabla\cdot(-\bv^{n+\demi} U)=0.
\label{eq:4}
\end{equation}
\item [ii)] ~Forward convection: 
\begin{equation}
\partial_t U+ \nabla\cdot(\bv^{n+\demi} U)=0.
\label{eq:5}
\end{equation}
\end{itemize}
Each convection problem is well-posed on the time interval $[0,\Delta t^n]$.
Let us now focus into these two steps and the way to solve them.
\subsection{Backward convection in Lagrangian description}
%--------------------------------------------------------
%
After the Lagrange step, if we solve the backward convection problem~\eqref{eq:4} 
over a time interval~$\Delta t^n$ using a Lagrangian description, we have
\begin{equation}
|K| (U_\ell)_K^{n,\star} =  |K^{n+1,L}| (U_\ell)_K^{n+1,L}. 
\label{eq:6}
\end{equation}
Actually, from the cell $K^{n+1,L}$ we go back to the original cell $K$ with
conservation of the conservative quantities. For $\ell=1$ (conservation of mass), 
we have
\[
|K|\, \rho_K^{n,\star} =  |K^{n+1,L}|\, \rho_K^{n+1,L}
\]
showing the variation of density by volume variation.
For $\ell=2,3,4$, it is easy to see that both velocity and specific total energy
are kept unchanged is this step:
\[
\bu^{n,\star} =  \bu^{n+1,L}, \quad 
E^{n,\star} =  E^{n+1,L}.
\]
Thus, this step is clearly computationally inexpensive.
\subsection{Forward convection in Eulerian description}
%-----------------------------------------------------
%
From the discrete field $(U_K^{n,\star})_K$ defined on the Eulerian cells $K$, 
we then solve the forward convection  problem~\label{eq:5} over a time step~$\Delta t^n$
under an Eulerian description. A standard Finite Volume discretization of the problem
will lead to the classical time advance scheme
\begin{equation}
U_K^{n+1} = U^{n,\star}_K - \frac{\Delta t^n}{|K|}\
\sum_{A\subset\partial_K} |A|\, U_A^{n+\demi,\star}\, (\bv_A^{n+\demi}\cdot \nu_A)
\label{eq:7}
\end{equation}
for some interface values $U_A^{n+\demi,\star}$ defined from the local neighbor values
$U^{n,\star}_K$. We finally get the expected Eulerian values $U_K^{n+1}$ at time $t^{n+1}$. \\[2ex]
Notice that from~\eqref{eq:6} and~\eqref{eq:7} we have also
\begin{equation}
|K|\,U_K^{n+1} = |K^{n+1,L}|\, U_K^{n+1,L} - \Delta t^n\
\sum_{A\subset\partial K} |A|\, U_A^{n+\demi,\star}\, (\bv_A^{n+\demi}\cdot \nu_A)
\label{eq:8}
\end{equation}
thus completely defining the remap step under the finite volume scheme form~\eqref{eq:8}.
We find that we no more need any mesh intersection or geometric consideration to 
achieve the remapping process. The finite volume form~\eqref{eq:8} is now suitable
for a straightforward vectorized SIMD treatment. From~\eqref{eq:8} it is easy to achieve second-order accuracy for the remapping step by usual finite volume tools
(MUSCL reconstruction + second-order accurate time advance scheme for example).
\subsection{Full Lagrange+remap time advance}
%---------------------------------------------
%
Let us note that the Lagrange+remap scheme is actually a conservative finite
volume scheme: putting~\eqref{eq:2} into~\eqref{eq:8} gives for all $\ell$:
\small
\begin{eqnarray}
(U_\ell)_K^{n+1} = (U_\ell)_K^n &-& \frac{\Delta t^n}{|K|}
\sum_{A^{n+\demi,L}\subset \partial K^{n+\demi,L}} |A^{n+\demi,L}|\, (\bpi_\ell)_A^{n+\demi,L}\cdot \nu_A^{n+\demi,L} 	\nonumber \\
&-& \frac{\Delta t^n}{|K|}\
\sum_{A\subset\partial K} |A|\, (U_\ell)_A^{n+\demi,\star}\, (\bv_A^{n+\demi}\cdot \nu_A)
\label{eq:9}
\end{eqnarray}
\normalsize
that can also be written
\small
\begin{eqnarray}
(U_\ell)_K^{n+1} = (U_\ell)_K^n &-& \frac{\Delta t^n}{|K|}
\sum_{A\subset\partial K}  |A|\,\left(\frac{|A^{n+\demi,L}|}{|A|}\, (\bpi_\ell)_A^{n+\demi,L}\cdot \nu_A^{n+\demi,L}\right) 	\nonumber \\
&-& \frac{\Delta t^n}{|K|}\
\sum_{A\subset\partial K}|A|\, \left( (U_\ell)_A^{n+\demi,\star}\, (\bv_A^{n+\demi}\cdot \nu_A)\right).
\label{eq:10}
\end{eqnarray}
\normalsize
We recognize in~\eqref{eq:10} pressure-related fluxes and convective fluxes that define
the whole numerical flux.
\subsection{Comments}\label{sec:44}
%----------------------------
%
The finite volume formulation~\eqref{eq:10} is attractive and seems rather
simple at first sight. But we should not forget that we have to compute a velocity
Lagrange vector field $\bv^{n+\demi}=\bu^{n+\demi,L}$ where the variables
should be located at cell nodes to return a well-posed deformation.
Moreover, expression~\eqref{eq:10} involves geometric elements like
the length of the deformed edges~$A^{n+\demi,L}$.
Among the rigorous collocated Lagrangian solvers, let us mention the GLACE
scheme by Despr\'es-Mazeran \cite{Despres2005} and the cell-centered EUCCLHYD
solver by Maire et al~\cite{Maire2007}. Both are rather computationally
expensive and their second-order accurate extension is not easy to achieve.

Although it is possible to couple these Lagrangian solvers with the flux-balanced
remapping formulation, it is also of interest to think about ways to 
simplify or approximate the Lagrange step without losing second-order accuracy.
One of the difficulty in the analysis of Lagrange-remap schemes is that, in some sense,
space and time are coupled by the deformation process.

Below, we derive a formulation that leads to a clear separation between space and time, 
in order to simply control the order of accuracy. The idea is to make
the time step tend to zero in the Lagrange-remap scheme (method of lines \cite{Schiesser1991}), 
then exhibit the
instantaneous spatial numerical fluxes through the Eulerian cell edges 
that will serve for the construction of an explicit finite volume scheme. Because the method needs an approximate Riemann solver in Lagrangian form, we will call it
a Lagrange-flux scheme.
\section{Derivation of a second-order accurate Lagrange-flux scheme}\label{sec:4}
%------------------------------------------------------------------
%
From the intermediate conclusions of the discussion~\ref{sec:44} above, we would like
to be free from any rather expensive collocated Lagrangian solver. However, such
a Lagrangian solver seems necessary 
to correctly and accurately define the deformation velocity field $\bv^{n+\demi}$ at time
$t^{n+\demi}$. \medskip

In what follows, we are trying to deal with time accuracy in a different manner. 
Let us come back to the Lagrange+remap formula~\eqref{eq:10}. Let us consider
a ``small'' time step $t$ that fulfils the usual stability CFL condition.
We have
\small
\begin{eqnarray*}
(U_\ell)_K(t) = (U_\ell)_K^n &-& \frac{t}{|K|}
\sum_{A\subset\partial K}  |A|\,\left(\frac{|A^L(t/2)|}{|A|}\, (\bpi_\ell)_A^{L}(t/2)\cdot \nu_A^{L}(t/2)\right) 	\nonumber \\
&-& \frac{t}{|K|}\
\sum_{A\subset\partial K}|A|\, (U_\ell)_A^{\star}(t/2)\, \bv_A(t/2)\cdot \nu_A.
\end{eqnarray*}
\normalsize
By making $t$ tend to zero, ($t>0$), we have $A^L(t/2)\rightarrow A$,
$(\bpi_\ell)^L(t/2)\rightarrow \bpi_\ell$,
$\bv(t/2)\rightarrow\bu$, $(U_\ell)^\star \rightarrow U_\ell$,
then we get a semi-discretization in space of the conservation laws. That can be seen as
a particular method of lines (\cite{Schiesser1991}):
\small
\begin{equation}
\frac{d (U_\ell)_K}{dt} = - \frac{1}{|K|}
\sum_{A\subset\partial K}  |A|\, ((\bpi_\ell)_A\cdot \nu_A) 
-\frac{1}{|K|}\
\sum_{A\subset\partial K}|A|\, (U_\ell)_A\, (\bu_A\cdot\nu_A).
\label{eq:11}
\end{equation}
\normalsize
We get a classical finite volume method
\[
\frac{dU_K}{dt} = - \frac{1}{|K|}\,\sum_{A\subset \partial K} |A| \ \Phi_A
\]
with a numerical flux $\Phi_A$ whose components are
\begin{equation}
(\Phi_\ell)_A =  (U_\ell)_A\, (\bu_A\cdot\nu_A) + (\bpi_\ell)_A\cdot \nu_A. 
\label{eq:12}
\end{equation}
In~\eqref{eq:11}, pressure fluxes $(p_\ell)_A$ and interface 
normal velocities $(\bu_A\cdot \nu_A)$
can be computed from an approximate Riemann solver in Lagrangian coordinates
(for example the Lagrangian HLL solver, see~\cite{Toro}). Then, the interface
states $(U_\ell)_A$ should be computed from an upwind process according to the sign
of the normal velocity $(\bu_A\cdot\nu_A)$. This is interesting because the resulting flux has similarities with the so-called advection upstream
splitting method (AUSM) flux family proposed by Liou~\cite{Liou}, but the construction here is different and, in some sense, justifies the AUSM 
splitting. \medskip

To get higher-order accuracy in space, one can
use a standard MUSCL reconstruction + slope limiting process involving 
classical slope limiters like for example Sweby's limiter 
function~\cite{Sweby1984}: 
\begin{equation}
\phi(a,b) = (ab>0) \ \text{sign}(a) \ \max\big( 
\min(|a|,\beta |b|),\ \min(\beta |a|, |b|)\big)
\label{eq:sweby}
\end{equation}
with $\beta\in [1,2]$ for achieving second order accuracy.
At this stage,
because there is no time discretization, everything is defined on the Eulerian mesh and fluxes
are located at the edges of the Eulerian cells. This is one originality of this scheme compared to the legacy staggered Lagrange-remap scheme that has to use variables defined on the Lagrangian cells. \medskip

To get high-order accuracy in time, one can then apply a 
standard high-order time advance scheme (RK2, etc.). For the second-order Heun scheme
for example, we have the following algorithm:
\begin{enumerate}
\item Compute the time step $\Delta t^n$ subject to some CFL condition;
\item \textbf{Predictor step}. MUSCL reconstruction + slope limitation on primitive variables $\rho$, $\bu$ and $p$. From the discrete states~$U_K^n$,
compute a discrete gradient for each cell $K$.
\item Use a Lagrangian approximate Riemann solver to compute pressure fluxes $\bpi_A^n$ and interface velocities $\bu_A^n$
\item Compute the upwind edge values $(U_\ell)_A^n$ according to the sign of
$(\bu_A^n\cdot \nu_A)$;
\item Compute the numerical flux $\Phi_A^n$ as defined in~\eqref{eq:12};
\item Compute the first order predicted states $U_{K}^{\star,n+1}$:
\[
U_K^{\star,n+1} = U_K^n - \frac{\Delta t^n}{|K|}\sum_{A\subset \partial K}
|A| \ \Phi_A^n
\]
\item \textbf{Corrector step}. MUSCL reconstruction + slope limitation: from the discrete values~$U_K^{\star,n+1}$,
compute a discrete gradient for each cell $K$.
\item Use a Lagrangian approximate Riemann solver to compute pressure fluxes $\bpi_A^{\star,n+1}$ and interface velocities $\bu_A^{\star,n+1}$
\item Compute the upwind edge values $(U_\ell)_A^{\star,n+1}$
according to the sign of $(\bu_A^{\star,n+1}\cdot \nu_A)$;
\item Compute the numerical flux $\Phi_A^{\star,n+1}$ as defined in~\eqref{eq:12};
\item Compute the second-order accurate states $U_{K}^{n+1}$ at time $t^{n+1}$:
\[
U_K^{n+1} = U_K^n - \frac{\Delta t^n}{|K|}\sum_{A\subset \partial K}
|A|\ \frac{\Phi_A^n + \Phi_A^{\star,n+1}}{2}.
\]
\end{enumerate}
One can appreciate the simplicity of the numerical solver compared to the legacy staggered Lagrange-remap algorithm. The complexity of the latter mainly due to various kernel (function) calls and too much communications is detailed in~\cite{PARCO2015}. Here the predictor and corrector kernel functions have similar
programming codes and there is no intermediate variables to save in memory.
\subsection{Lagrangian HLL approximate solver}
%--------------------------------------------
%
A HLL approximate Riemann solver \cite{Toro} in Lagrangian coordinates
can be used to easily
compute interface pressure and velocity. For a local Riemann problem made
of a left state~$U_L$ and a right state $U_R$, the contact pressure $p^\star$
is given by the formula
\begin{equation}
p^\star = \frac{\rho_R p_L + \rho_L p_R}{\rho_L+\rho_R}
- \frac{\rho_L\rho_R}{\rho_L+\rho_R}\, \max(c_L,\,c_R)\, (u_R-u_L),
\label{eq:13}
\end{equation}
and the normal contact velocity $u^\star$ by
\begin{equation}
u^\star =  \frac{\rho_L u_L + \rho_R u_R}{\rho_L+\rho_R}
- \frac{1}{\rho_L+\rho_R}\, \frac{p_R-p_L}{\max(c_L,\,c_R)}
\label{eq:14}
\end{equation}
leading to simple formulas easily implementable in the second-order Heun time integration scheme.
\subsection{Numerical experiments}
%--------------------------------
%
As example, we test the Lagrange-flux scheme presented in section~\ref{sec:4} 
on few one-dimensional shock tube problems. We use a Runge-Kutta 2 (RK2) time integrator and a MUSCL reconstruction with the Sweby slope 
limiter given in~\eqref{eq:sweby}.

\begin{figure*}
%\notoprule         % Very first command!
%!\centering\includegraphics[height=0.16\textheight]{sod384.png}
%\centering
\hspace{-0.11\textwidth}\includegraphics[width=1.2\textwidth]{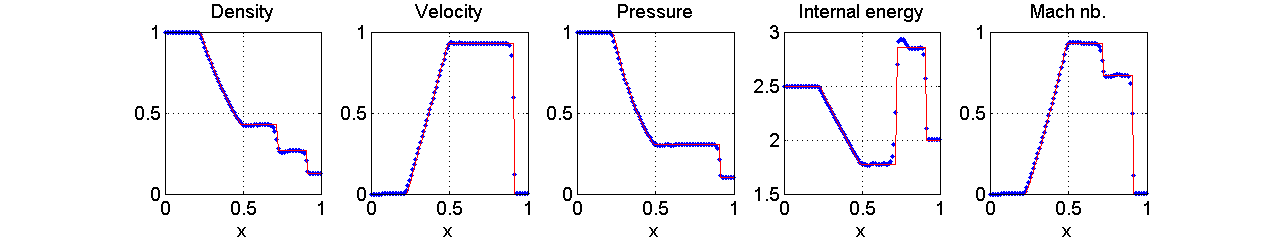}

\hspace{-0.11\textwidth}\includegraphics[width=1.2\textwidth]{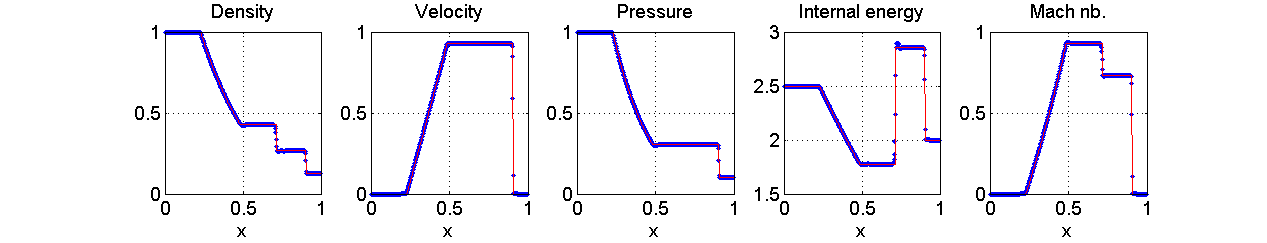}
\caption{Second-order Lagrange-flux scheme on reference Sod's hock tube problem
on two different mesh grids, 100 and 400 respectively. The solid red line is the analytic solution and blue points are the numerical solution.} \label{fig:sod}
\end{figure*}

\begin{figure*}
%\notoprule         % Very first command!
%\begin{center}
\hspace{-0.11\textwidth}\includegraphics[width=1.2\textwidth]{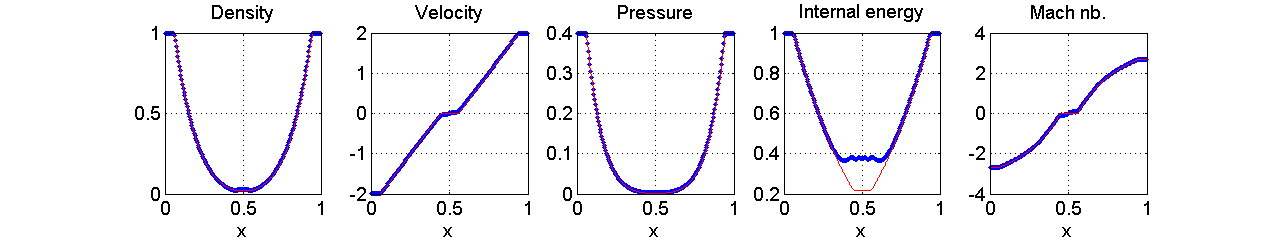}

\hspace{-0.11\textwidth}\includegraphics[width=1.2\textwidth]{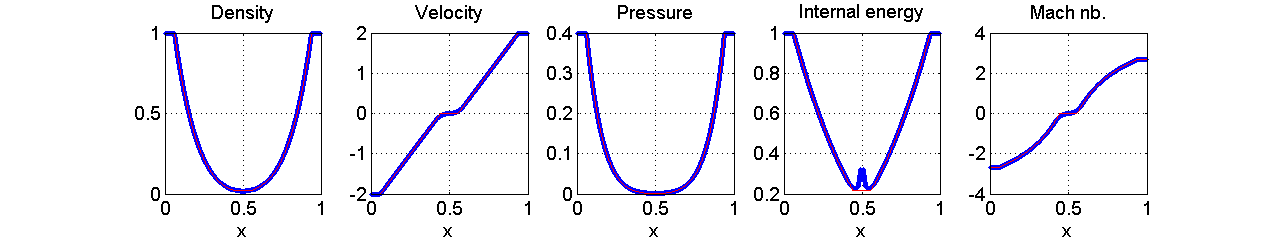}
\caption{Second-order Lagrange-flux scheme on a double rarefaction near-vacuum case on two different mesh grids, 200 and 2000 points respectively. The solid red line is the analytic solution.} \label{fig:cavit}
%\end{center}
\end{figure*}

\begin{figure*}
\hspace{-0.11\textwidth}\includegraphics[width=1.2\textwidth]{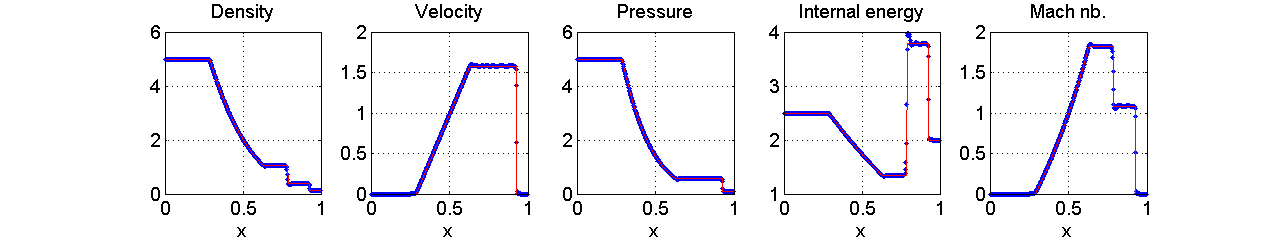}
\caption{Shock tube problem with sonic rarefaction fan, 400 mesh points, final time $T$ is $T=0.18$, the solid red line is the analytic solution.}
\label{fig:4}
\end{figure*}

\begin{figure*}
\hspace{-0.11\textwidth}\includegraphics[width=1.2\textwidth]{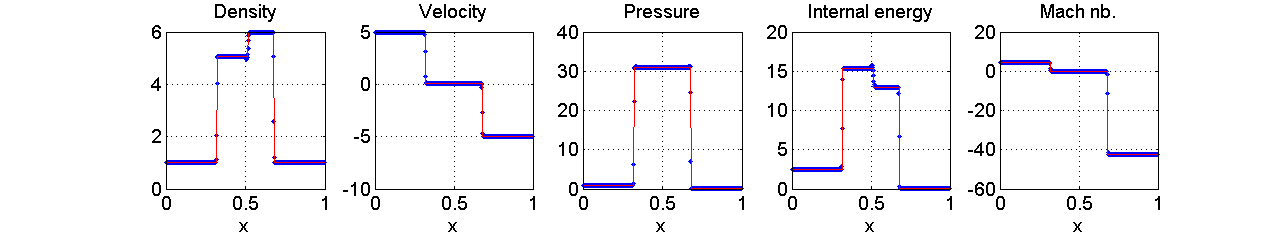}
\caption{Shock-shock case with subsonic-hypersonic shock, 400 mesh points,  final 
time $T$ is $T=0.16$, the solid red line is the analytic solution.}
\label{fig:5}
\end{figure*}

\paragraph{Sod's shock tube \cite{Sod71}} 
The initial data defined on space interval~$[0,1]$ is made of two constant
states $(\rho,u,p)_L = (1,0,1)$ and $(\rho,u,p)_R=(0.125,0,0.1)$ with initial discontinuity at $x=0.5$.
We successively test the method on two uniform mesh grids made of 100 and 400 cells respectively. 
The final computational time is $T=0.23$ and we use a CFL number equal to 0.25
and a limiter coefficient $\beta=1.5$.
On figure~\ref{fig:sod}, one can observe a nice behavior of the
Euler solver, with sharp discontinuities and a low numerical diffusion
into the rarefaction fan even for the coarse grid.
%
%\begin{figure*}
%\notoprule         % Very first command!
%\centering\includegraphics[height=0.16\textheight]{sod384.png}
%\caption{Second-order Lagrange-flux scheme on reference Sod's 1D shock tube %problem. Time is $T=0.23$, 384 mesh points. Here Sweby's slope limiter%
%with coefficient~1.5 is used.} \label{fig:sod}
%\end{figure*}
%
\paragraph{Two-rarefaction shock tube}
The second reference example is a case of two moving-away rarefaction fans under near-vacuum conditions (see Toro~\cite{Toro}). It is known that the Roe scheme breaks down for
this case. The related Riemann problem is made of the left state $(\rho,u,p)_L=(1,-2,0.4)$ and right state $(\rho,u,p)_R=(1,2,0.4)$. The final time of $T=0.16$. We again test the method of a coarse mesh (200 points) and a fine mesh (2000 points). Numerical results are given in 
figure~\ref{fig:cavit}. The numerical scheme appears to be robust especially in
near-vacuum zones where both density and pressure are close to zero.
%
%\begin{figure*}
%\notoprule         % Very first command!
%\begin{center}
%\includegraphics[width=\textwidth]{cavit.png}
%\caption{Second-order Lagrange-flux scheme on a double rarefaction near-vacuum %case. Final time is $T=0.16$, 8192 mesh points.} \label{fig:cavit}
%\end{center}
%\end{figure*}
%
\paragraph{Case with sonic rarefaction and supersonic contact}
The following shock tube case with initial data 
$(\rho,u,p)_L=(5,0,5)$ and $(\rho,u,p)_R=(0.125,0,0.1)$ generates
a sonic 1-rarefaction, a supersonic 2-contact discontinuity and a 3-shock wave.
The final time is $T=0.16$ and we use 400 mesh points, $CFL=0.25$. Numerical
results show a good capture of the rarefaction wave, without any
non-entropic expansion-shock (see figure~\ref{fig:4}). 
\paragraph{Case of shock-shock hypersonic shock tube}
This last shock tube problem is a violent flow case made of two hitting fluids
with $(\rho,u,p)_L=(1,5,1)$ and~$(\rho,u,p)_R=(1,-5,0.01)$. Both 1-wave and
3-wave are shock waves, and the right state has a Mach number of order 40.
Final time is $T=0.16$, we use 400 grid points and the limiter coefficient
$\beta$ is here~$1$ (equivalent to the minmod limiter). One can observe
of very nice behavior of the solver: there is no pressure or velocity
oscillations at the contact discontinuity, and the numerical scheme preserves
the positivity of density, pressure and internal energy (see figure~\ref{fig:5}). 
\section{Performance results}
%===========================
%
In this section, we compare the new Lagrange-flux scheme to the
reference (staggered) Lagrange-remap scheme in terms of millions of cell updates
per second (denoted hereafter as MCUPs). Tests are performed on a
standard $2\times 8$ cores Intel Sandy Bridge server E5-2670. Each core has a
frequency of 2.6~GHz, and supports Intel's AVX (Advanced Vector Extension) 
vector instructions. For multicore support, we use the multithreading
programming interface \texttt{OpenMP}. \medskip

In the reference staggered Lagrange-remap solver (see~\cite{PARCO2015}),
thermodynamic variables are defined at grid cell centers while
velocity variables are defined at mesh nodes. Due to this staggered
discretization and the alternating direction (AD) remapping
procedures, this solver is decomposed into nine kernels. This
decomposition mechanically decreases the mean arithmetic intensity (AI) of
the solver.

On the other hand, the Lagrange-flux algorithm consists in only two
kernels with a relative high arithmetic intensity which leads to two
compute-bound (CB) kernels. In the first kernel, named
\texttt{PredictionLagrangeFlux()}, an appropriate Riemann solver is called, face fluxes are
computed and variables are updated for the prediction. The second
kernel, named \texttt{CorrectionLagrangeFlux()}, is close in terms of
algorithmic steps, since it also uses a Riemann solver, computes fluxes
and updates the variables for the correction part of the solver.

In order to assess the scalability and absolute performance of both
schemes, we present in table~\ref{tab_perfs} a performance comparison
study.  First, we notice that the baseline performance --- e.g.\ the
single core absolute performance without vectorization --- is quite
similar for the two schemes, as can be seen in the first
column. However, we the Lagrange-flux scheme has a better scalability,
due to both vectorization and multithreading: our Lagrange-flux
implementation achieves a speedup of 31.1X with 16 cores and AVX
vectorization (while ideal speed-up is 64) whereas the reference
Lagrange-remap algorithm reaches a speed-up of only 14.8X. This
difference is mainly due to the memory-bound kernels composing the
reference Lagrange-remap scheme. Indeed, speedups due to AVX
vectorization and multithreading are not ideal for kernels with
relatively low intensity since memory bandwidth is shared between cores.
\begin{table}
  \begin{center}
  \label{table:multicore}
  \begin{tabular}{@{}l c c c c @{}}
  \hline
 \textbf{Scheme}& 1 core & 1 core & 16 cores & Scalability\\
 &        & AVX & AVX \\[2ex]
%\hline
Lagrange-flux & 2.6& 5.8 &  81.0 & 31.1\\
%\hline
Reference & 2.5& 3.8 & 37.0& 14.8\\
% %\hline
\hline
\end{tabular}
  \caption{Performance comparison between the reference Lagrange-remap solver and
      the Lagrange-flux solver in millions of cell updates per second
      (MCUPs), using different machine configurations. Scalability (last
      column) is computed as the speedup of the multithreaded
      vectorized version compared to the baseline purely
      sequential version. Tests are performed for fine meshes, such
      that kernel data lies in DRAM memory. The Lagrange-flux solver
      exhibits superior scalability, because it has --- by design ---
      better arithmetic intensity.}\label{tab_perfs}
\end{center}
\end{table}
\section{Dealing with multimaterial flows}\label{sec:5}
%----------------------------------------
%
Although this is not the aim and the scope of the present paper, we would like to give
an outline of the possible extension of Lagrange-flux schemes to compressible
multimaterial/multifluid flows, i.e. flows that are composed of different immiscible
fluids and separated by free boundaries. \\
For pure Lagrange+remap schemes, usually VOF-based interface reconstruction (IR) algorithms
are used (Young's PLIC, etc.). After the Lagrangian evolution, for the cells that
host more than one fluid, fluid interfaces are reconstructed. 
During the remapping step, one has to evaluate the mass fluxes per material.
From the computational point of view and computing performance, this process generally slows down the whole performance because of many array indirections in memory and
specific treatment into mixed cells along with the material interfaces. \medskip

If the geometry of the Lagrangian cells is not completely known (as in the case of Lagrange-flux
schemes), anyway we have to proceed differently. A possibility is to use interface
capturing (IC) schemes, e.g.\ conservative Eulerian schemes that evaluate the convected mass fluxes through Eulerian cell edges. This can be achieved by the
use of antidiffusive/low-diffusive advection solvers in the spirit of 
Despr\'es-Lagouti\`ere's limited-downwind scheme~\cite{Despres2001} of 
VoFire~\cite{Despres2010}. In a recent work~\cite{DeVuyst2015}, we have analyzed the origin
of known artifacts and numerical interface instabilities for this type of solvers and concluded that the reconstruction of fully multidimensional gradients with multidimensional gradient limiters was  necessary. Thus, we decided to use low-diffusive advection schemes with a Multidimensional Limiting Process (MLP) in the spirit of
\cite{Park2011}. The resulting method is quite accurate, shape-preserving and
free from any artifact. We show some numerical results in the two next subsections.
Let us emphasize that the interface capturing strategy perfectly fits with the Lagrange-flux flow description, and the resulting schemes are really suitable for vectorization
(SIMD feature) with data coalescence into memory.
\subsection{Interface capturing for pure advection problems}\label{sec:rider}
\begin{figure}[h!]
		\begin{subfigure}{0.5\textwidth}
		\includegraphics[width=0.85\linewidth]{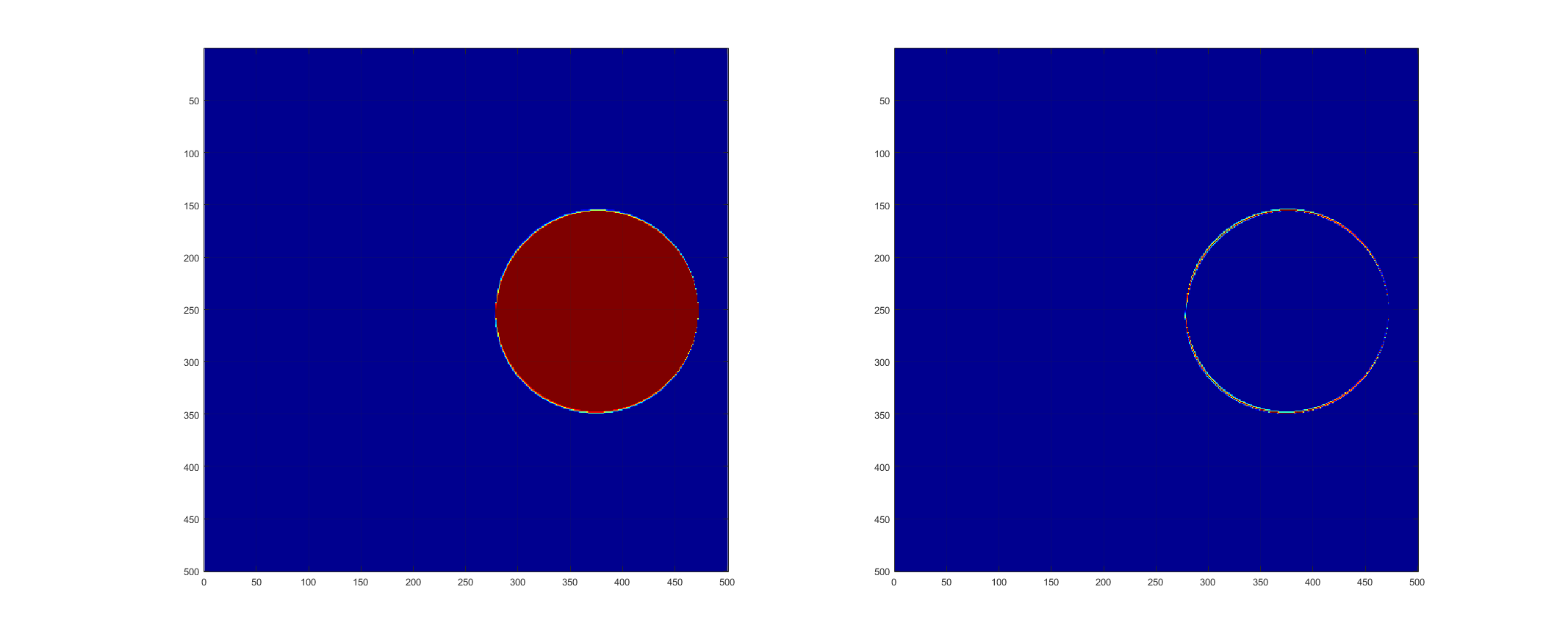}
		\caption{At initial time $t=0$}
		\end{subfigure}
		\begin{subfigure}{0.5\textwidth}
		\includegraphics[width=0.85\linewidth]{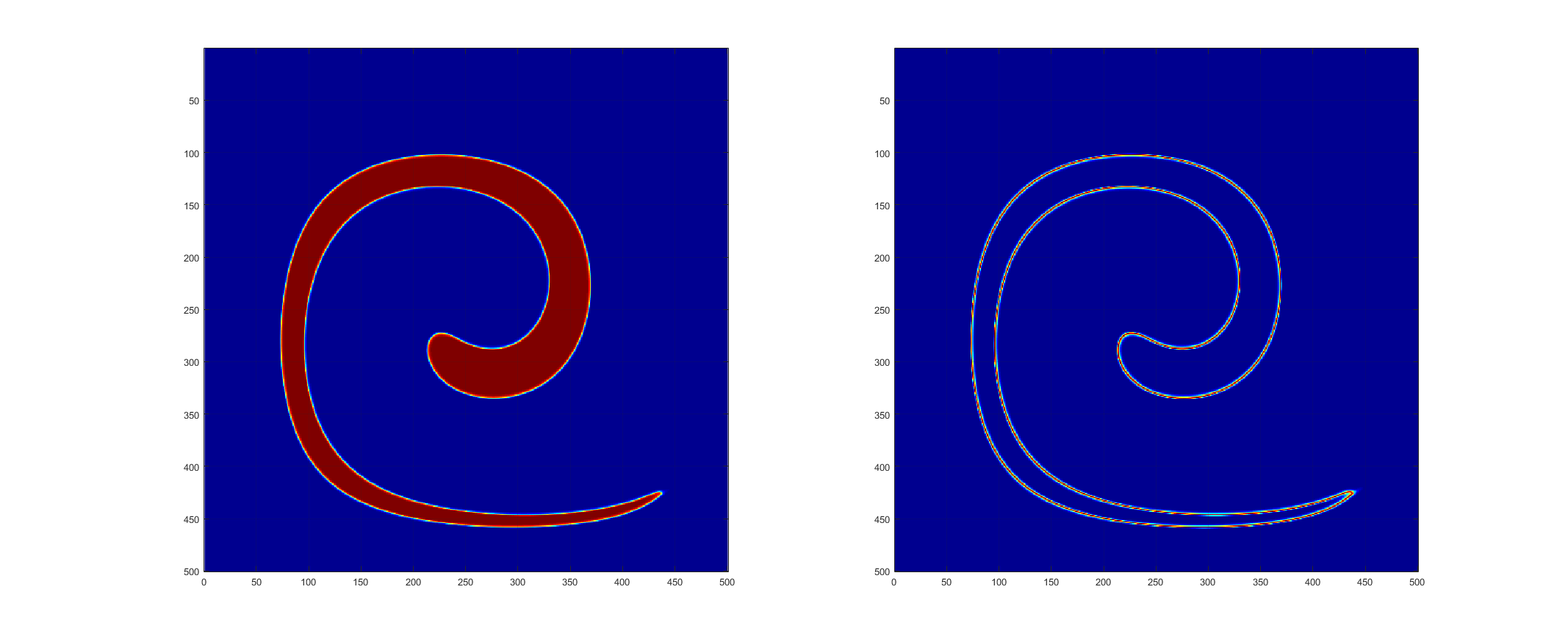}
		\caption{At time $t=3$}
		\end{subfigure}
		\begin{subfigure}{0.5\textwidth}
		\includegraphics[width=0.85\linewidth]{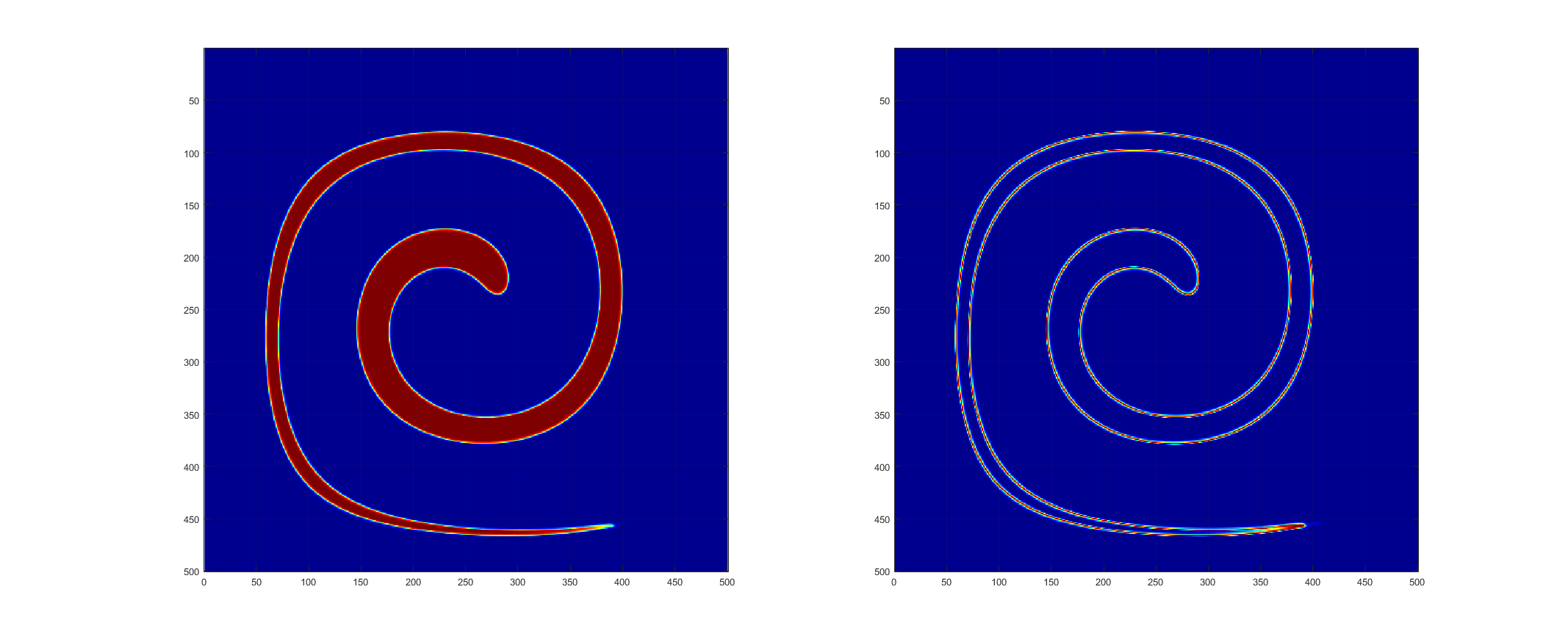}
		\caption{At time $t=6$}
		\end{subfigure}
		\begin{subfigure}{0.5\textwidth}
		\includegraphics[width=0.85\linewidth]{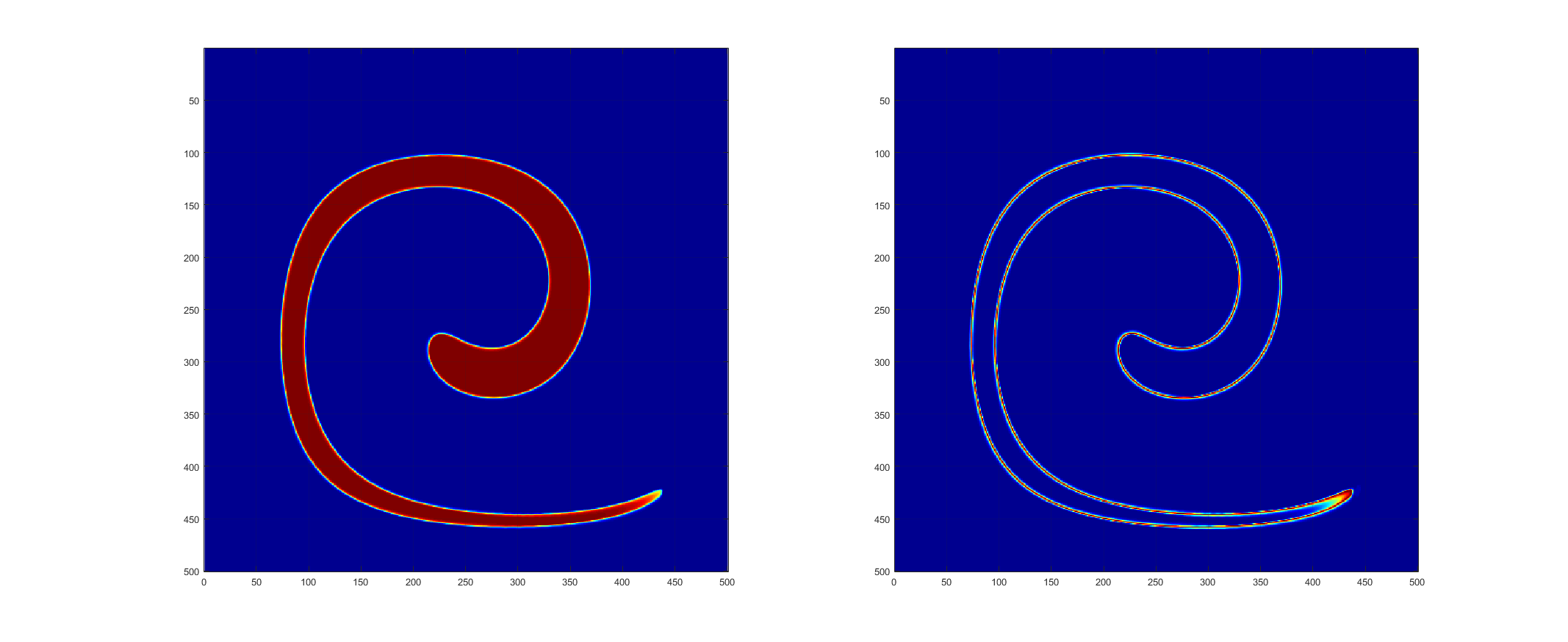}
		\caption{At time $t=9$}
		\end{subfigure}
		\begin{subfigure}{0.5\textwidth}
		\includegraphics[width=0.85\linewidth]{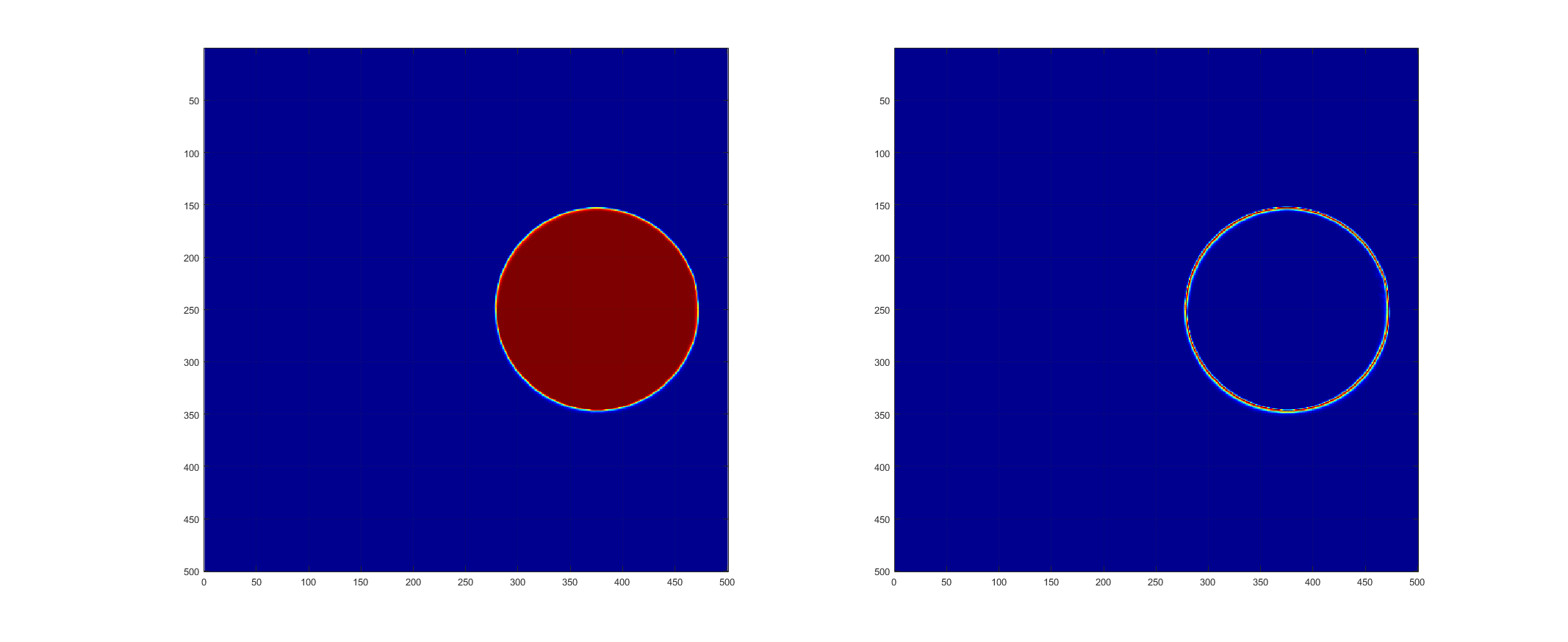}
		\caption{At final time $t=12$}
		\end{subfigure}
		\begin{subfigure}{0.5\textwidth}
		\includegraphics[width=0.85\linewidth]{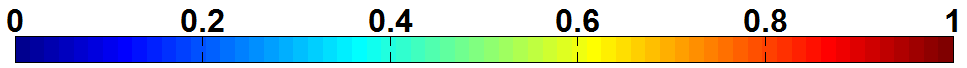}
		\end{subfigure}
		\caption{Validating the low-diffusive interface capturing scheme on the Kothe-Rider advection case, mesh $500^2$.} \label{fig:kr2}
		\end{figure}

%----------------------------------------------------------
%
Let us first present numerical results on a pure scalar linear advection problem.
The forward-backward advection case proposed by Rider and Kothe~\cite{Rider}
is a hat-shaped function which is advected and stretched into a rotating vector
field, leading to a filament structure. Then by applying the opposite velocity
field, one have to retrieve the initial disk shape. In figure~\ref{fig:kr2}, we 
show the numerical solutions obtained on a grid $500^2$ for both the passive scalar
field of variable $z\in [0,1]$ and the quantity $z(1-z)$ that indicates
the numerical spread rate of the diffuse interface. One can conclude the good behavior of the method, providing both stability, accuracy and shape preservation.
\subsection{Three-material hydrodynamics problem}\label{sec:triple}
%------------------------------------------------		
%
We then consider our interface capturing method for multifluid
hydrodynamics. Because material mass fractions are advected, that is
\[
\partial_t y_k + \bu\cdot\nabla y_k = 0, \quad k=1,...,N,
\]
one can use the advection solver of these variables but we prefer dealing with the conservative form of the equations
\[
\partial_t (\rho y_k) + \nabla\cdot(\rho y_k \bu) = 0
\]
in order to enforce mass conservation (see also \cite{Champmartin2014}). 
It is known that Eulerian interface-capturing schemes generally
produce spurious pressure oscillations at material interfaces (\cite{Abgrall,Saurel}).
Some authors propose locally non conservative approaches \cite{Farhat,Bachmann} to prevent
from any pressure oscillations. Here we have a full conservative Eulerian strategy
involving a specific limiting process which is free from any pressure
oscillation at interfaces, providing strong robustness. This will be explained
in a next paper.

\begin{figure}
%		\notoprule         % Very first command!
		\centering\includegraphics[width=0.9\linewidth]{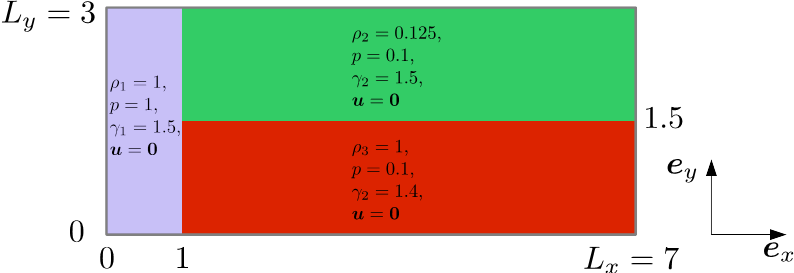}
		\caption{Geometry and initial configuration for the reference triple-point case.}
		\label{fig:triple}
		\end{figure}

\begin{figure}
		%\centering\includegraphics[width=0.99\textwidth]{Triple}
		\begin{subfigure}{0.49\textwidth}
		\begin{center}
		\includegraphics[width=0.9\linewidth]{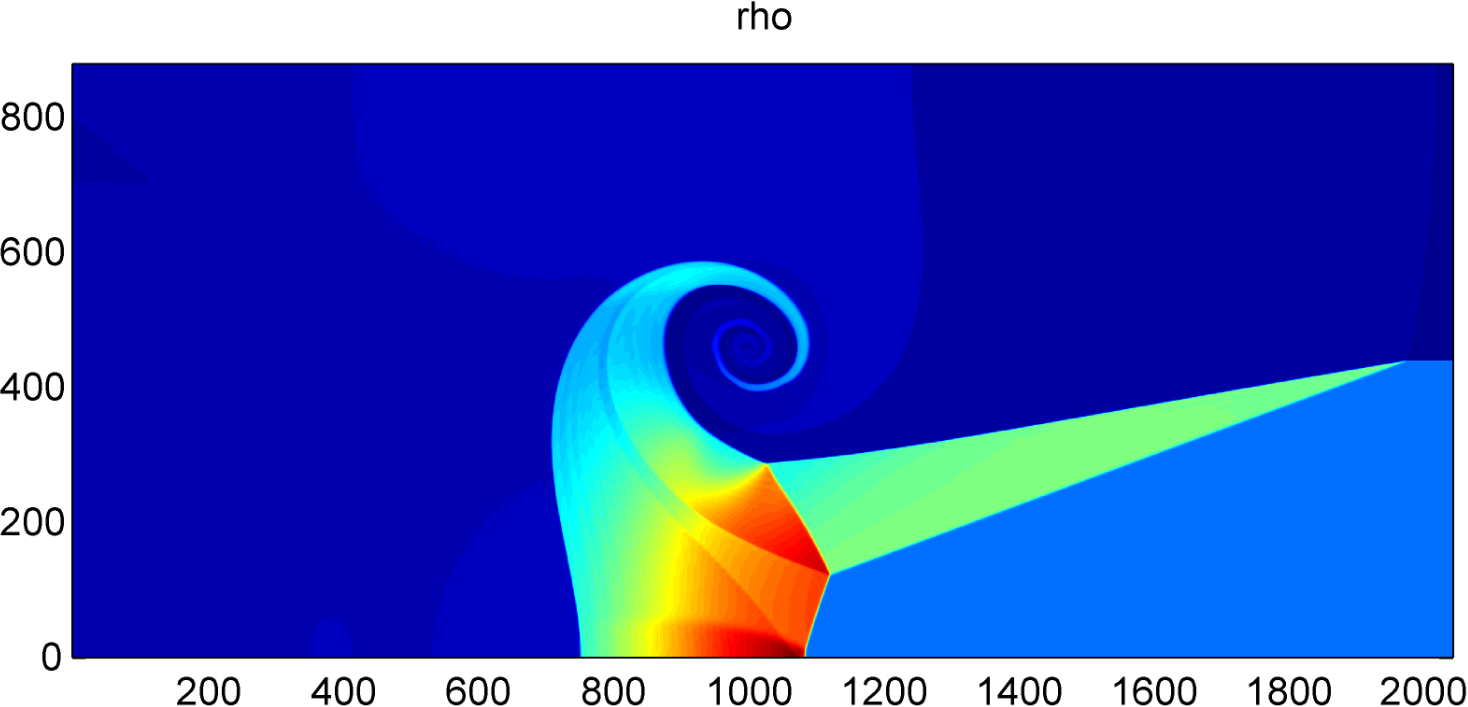}
		\end{center}
		\caption{Density field}
		\end{subfigure}
		\begin{subfigure}{0.49\textwidth}
		\begin{center}
		\includegraphics[width=0.9\linewidth]{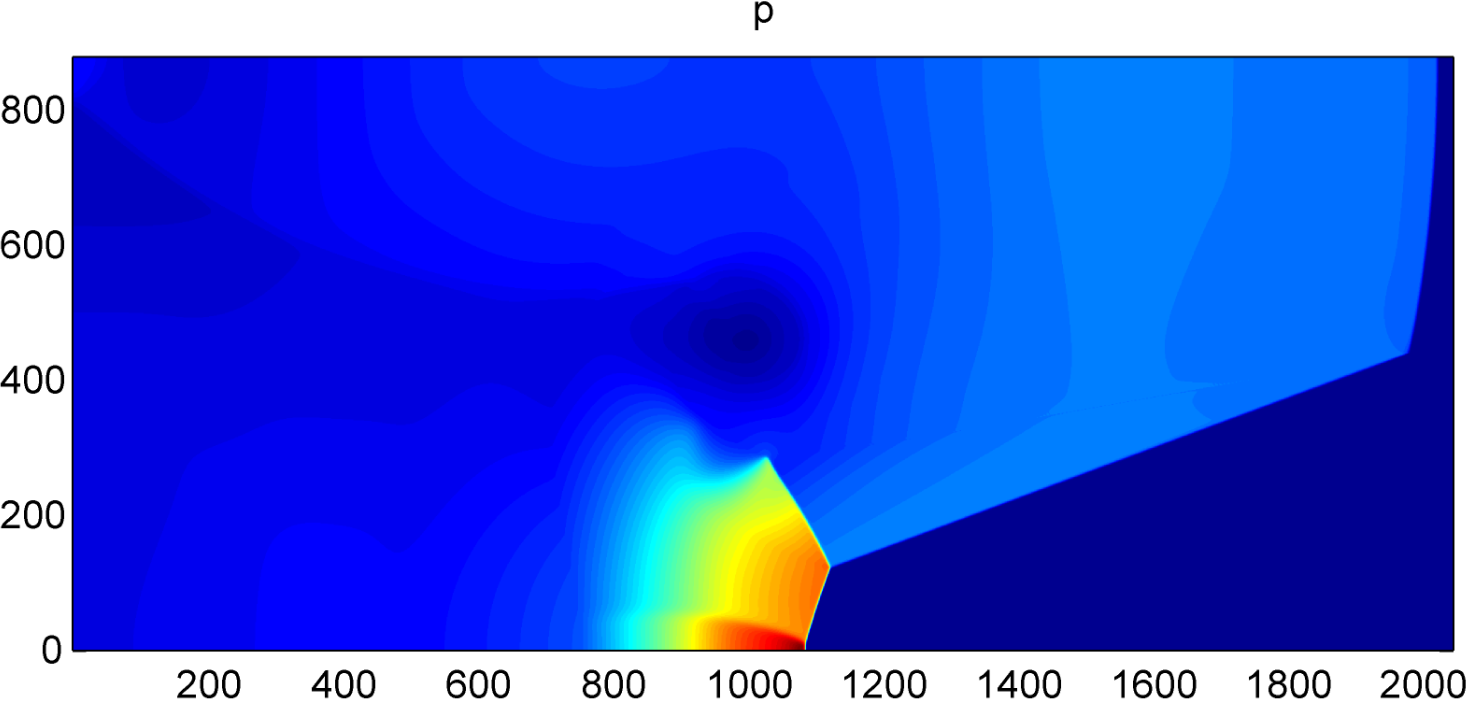}
		\end{center}
		\caption{Pressure field}
		\end{subfigure}
		\begin{subfigure}{0.49\textwidth}
		\begin{center}
		\includegraphics[width=0.9\linewidth]{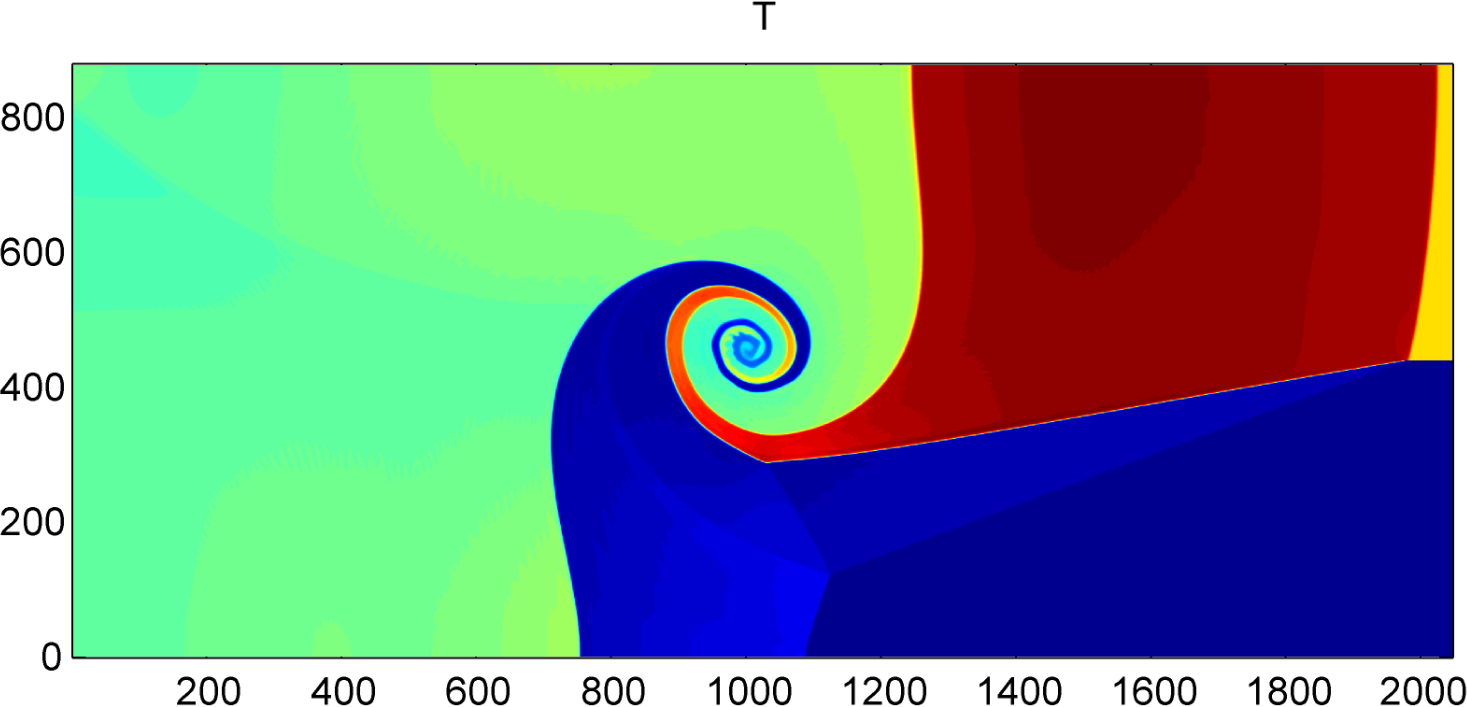}
		\end{center}
		\caption{Temperature field}
		\end{subfigure}
		\begin{subfigure}{0.49\textwidth}
		\begin{center}
		\includegraphics[width=0.9\linewidth]{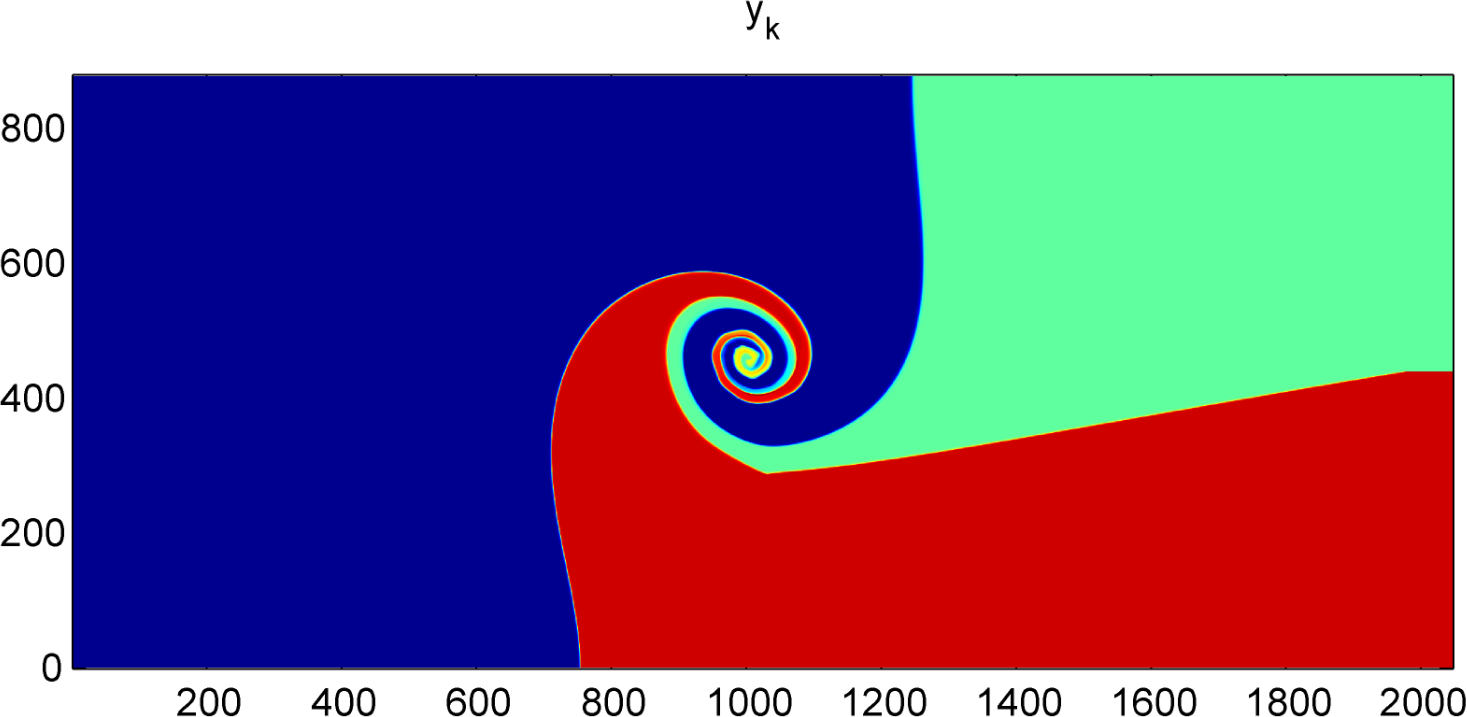}
		\end{center}
		\caption{Colored representation of material indicators}
		\end{subfigure}
		\caption{Results on the multimaterial ``triple point'' case  (perfect gases) 
		using a collocated Lagrange+remap solver 
		+ low-diffusive interface capturing advection scheme, mesh made of 
		$2048\times878$ points. Final time is $T=3.3530$.}
		\label{fig:triple1}
		\end{figure}

The multimaterial Lagrange-flux scheme is tested on the reference ``triple point''
test case, found e.g. in Loub\`ere et al.~\cite{Loubere}. 
This problem is a three-state two-material 2D Riemann problem in a rectangular vessel.
The simulation domain is $\Omega=(0,7)\times(0,3)$ as described in figure~\ref{fig:triple}.
The domain is splitted up into three regions $\Omega_i$, $i=1,2,3$ filled with two perfect gases leading to a two-material problem. Perfect gas equations of state are used with $\gamma_1=\gamma_3=1.5$ and~$\gamma_2=1.4$. Due to the density differences, two
shocks in sub-domains $\Omega_2$ and $\Omega_3$ propagate with different speeds. This
create a shear along the initial contact discontinuity and the formation of a vorticity.
Capturing the vorticity is of course the difficult part to compute. We use a rather
fine mesh made of $2048\times 878$ points (about 1.8M cells).
On figure~\ref{fig:triple1}, we plot the density, pressure, temperature fields
respectively and indicate the location of the three material zones. One can observe
a nice capture of both shocks and contact discontinuities. The vortex is
also captured accurately.
		% TRIPLE
				%
\section{Concluding remarks and perspectives}\label{sec:6}
%-------------------------------------------
%
This paper is primarily focused on the redesign of Lagrange-remap hydrodynamics
solvers in order to achieve better HPC node-based performance. We have reformulated
the remapping step under a finite volume flux balance, allowing for a full SIMD
algorithm. As an unintended outcome, the analysis has lead us to the discovery of a new
promising family of Eulerian solvers -- the so-called Lagrange-flux solvers -- that show
simplicity of implementation, accuracy, and flexibility with a high-performance capability compared to the legacy staggered Lagrange-remap scheme.
Interface capturing methods can be easily plugged for solving multimaterial 
flow problems. Ongoing work is focused of the effective performance modeling, analysis
and measurement of Lagrange-flux schemes with comparison of reference ``legacy''
Lagrange-remap solvers including multimaterial interface capturing on different multicore processor architectures. Because of the multicore+vectorization scalability of Lagrange-flux schemes, one can also expect high-performance on manycore co-processors like
graphics processing units (GPU) or Intel MIC. This will be the aim of  next
developments.  
%
% 1. We do not claim that
% 2. Tried to improve the whole complexity in terms of FP operations + smaller
% of data to transfer -> Refactoring/rethinking
% 3. Satisfying the requirements also lead to the design of new numerical methods
% -- the so-called Lagrange-flux schemes
% 4. Extension to MM okay
%
%==============================================

%

%==============================================================================================================

\begin{thebibliography}{}
%
\bibitem{PARCO2015} R. Poncet, M. Peybernes, T. Gasc and F. De Vuyst,
Performance modeling of a compressible hydrodnamics solver on multicore CPUs, 
Proceedings of the int. conf. on parallel computing PARCO2015, Edinburgh, 2015 (in press).
%
\bibitem{Hirt1974} C.W. Hirt, A.A. Amsden and J.L. Cook, An arbitrary Lagrangian–Eulerian computing method for all flow speeds. Journal of Computational Physics, 14:227--253 (1974).
%
\bibitem{benson} D. Benson, Computational methods in Lagrangian and Eulerian
hydrocodes, CMAME, 99(2-3), 235--394 (1992).
%
\bibitem{youngs} D. Youngs, The Lagrange-remap method, in \textit{Implicit
Large Eddy Simulation: computing turbulent flow dynamics}, F.F. Grinstein,
L.G. Margolin and W.J. Rider (eds), Cambridge University Press (2007).
%
\bibitem{Williams2009} S. Williams, A. Waterman, and D. Patterson: Roofline: An Insightful Visual Performance Model for Multicore Architectures, Commun. ACM, 52, pp 65–76 (2009).
%
\bibitem{Treibig2010} J. Treibig and G. Hager, Introducing a Performance Model for Bandwidth-Limited Loop Kernels. Proceedings of the Workshop “Memory issues on Multi- and Manycore Platforms” at PPAM 2009, Lecture Notes in Computer Science, 6067, pp. 615--624 (2010).
%
\bibitem{Stengel2015} H. Stengel, J. Treibig, G. Hager and G. Wellein, Quantifying performance bottlenecks of stencil computations using the Execution-Cache-Memory model. Proc. ICS15, the 29th Int. Conf. on Supercomputing,  2015, DOI: 10.1145/2751205.2751240.
%
\bibitem{Collela1984} P. Colella and P.R. Woodward, The numerical simulation of two-dimensional fluid flow with strong shocks, J. Comput. Phys.,54:115--173 (1984).
%
\bibitem{Despres2005} B. Despr\'es and C. Mazeran, Lagrangian gas dynamics in two dimensions and Lagrangian systems, Arch. Rational Mech. Anal. 178 (2005) 327--372.
%
\bibitem{Maire2007} P.-H. Maire, R. Abgrall, J. Breil and J. Ovadia, A cell-centered Lagrangian scheme for compressible flow problems, SIAM J. Sci. Comput. 29 (4) (2007) 
1781--1824.
%
\bibitem{Maire2009} P.-H. Maire, A high-order cell-centered Lagrangian scheme for two-dimensional compressible fluid flows on unstructured meshes, J. Comput. Phys. 228
(2009) 2391--2425.
%
\bibitem{Dukowicz2000} J. K. Dukowicz and J. R. Baumgardner, Incremental
remapping as a transport/advection algorithm. J. Comput. Phys., 160, 318–335 (2000).
%
\bibitem{Schiesser1991} W. E. Schiesser, The Numerical Method of Lines, Academic Press, ISBN 0-12-624130-9 (1991).
%
\bibitem{Toro} E.F. Toro, Riemann solvers and numerical methods for fluid dynamics, 3rd Edition, Springer (2010).
%
\bibitem{Liou} M.S. Liou, A sequel to AUSM: AUSM+, J. Comp. Phys., 129(2),
364--382 (1996).
%
\bibitem{Sod71} G.A. Sod, A Survey of Several Finite Difference Methods for Systems of Nonlinear Hyperbolic Conservation Laws" . J. Comput. Phys. 27: 1--31 (1971).
%
\bibitem{Sweby1984} P.K. Sweby, High resolution schemes using flux-limiters for hyperbolic conservation laws, SIAM J. Num. Anal. 21 (5): 995--1011 (1984).
%
\bibitem{Despres2001} B. Despr\'es and F. Lagouti\`ere, Contact discontinuity capturing schemes for linear advection and compressible gas dynamics, J. Sci. Comp., 16(4), 479--524 (2001).
%
\bibitem{Despres2010} B. Despr\'es, F. Lagouti\`ere, E. Labourasse and I. Marmajou, An antidissipative transport scheme on unstructured meshes for multicomponent flows, IJFV, 30--65 (2010).
%
\bibitem{DeVuyst2015} F. De Vuyst, M. B\'echereau, T. Gasc, R. Motte, M. Peybernes
and R. Poncet, Stable and accurate low-diffusive  interface capturing advection schemes, submitted to proc. of the MULTIMAT2015 Conference W\"ursburg, special issue of the IJNMF (2015).
%
%
\bibitem{Park2011} J.S. Park and C. Kim, Multi-dimensional Limiting Process for Discontinuous Galerkin Methods on Unstructured Grids, chapter in Computational Fluid Dynamics 2010, Springer, 179--184 (2011).
%
\bibitem{Rider} A.J. Rider and D.B. Kothe, Reconstructing volume tracking, J. Comput. Phys., 141(2), 112--152 (1998).
%
\bibitem{Champmartin2014} A. Bernard-Champmartin and F. De Vuyst, A low diffusive Lagrange-remap scheme for the simulation of violent air-water free-surface flows, J. of Comput. Physics, 274, 19--49 (2014).
%
\bibitem{Abgrall} R. Abgrall, How to prevent pressure oscillations in multicomponent
flow calculations: a quasi conservative approach, J. Comp. Phys., 125(1),
150--160 (1996).
%
\bibitem{Saurel} R. Saurel and R. Abgrall, A simple method for compressible
multifluid flows, SIAM J. Sci. Comput., 21(3), 1115--1145 (1999).
%
\bibitem{Farhat} C. Farhat, A. Rallu and S. Shankaran, A higher-order
generalized ghost fluid method for the poor for the three-dimensional
two-phase flow computation of underwater implosions, J. Comp. Phys, 
227(16), 7640-7700 (2008).
%
\bibitem{Bachmann} M. Bachmann, P. Helluy, J. Jung, H. Mathis and
S. M\"uller, Random sampling remap for compressible two-phase flows,
Computers and fluids, 86, 275--283 (2013).
%
\bibitem{Loubere} R. Loub\`ere, P.-H. Maire, M. Shashkov, J. Breil and S. Galera, ReALE: A reconnection-based arbitrary-Lagrangian-Eulerian method, J. Comp. Phys., 229, 4724--4761 (2010).
%
%
\end{thebibliography}
\end{document}